\documentclass[12pt]{amsart}
\usepackage[colorlinks=true, pdfstartview=FitV, linkcolor=blue, citecolor=green, urlcolor=black,filecolor=magenta]{hyperref}

\usepackage{graphicx, verbatim,amsmath,amssymb}

\newcommand{\aaa}{{\mathcal A}}

\newcommand{\lcal}{{\mathcal L}}

\newcommand{\R}{{\mathbb R}}

\newcommand{\Z}{{\mathbb Z}}

\newtheorem{thm}{Theorem}[section]
\newtheorem*{thm*}{Theorem}

\everymath{\displaystyle}

\theoremstyle{definition}
\newtheorem{ex}{Example}

\begin{document}
\title{Finitely balanced sequences and plasticity of 1-dimensional tilings}
\author{Lorenzo Sadun}
\date{\today}
%\version{5}

\address{Department of Mathematics\\The University of
Texas at Austin\\ Austin, TX 78712} \email{sadun@math.utexas.edu}
\thanks{This work is partially supported by NSF
grant  DMS-1101326} 
%\date{Jan 20, 2013}

\begin{abstract}
  We relate a balancing property of letters for bi-infinite sequences
  to the invariance of the resulting 1-dimensional tiling dynamics
  under changes in the lengths of the tiles. If the language of the
  sequence space is finitely balanced, then all length changes in the
  corresponding tiling space result in topological conjugacies, up to
  an overall rescaling.
\end{abstract}

\maketitle

\setlength{\baselineskip}{.6cm}

%\setcounter{tocdepth}{3}
%\makeatletter
%\def\l@subsection{\@tocline{2}{0pt}{25pt}{5pc}{}}
%\def\l@subsubsection{\@tocline{2}{0pt}{50pt}{5pc}{}}
%\makeatother
%\tableofcontents

\section{Introduction and statement of result}

Spaces of bi-infinite sequences are closely related to spaces of
1-dimensional tilings. Given a tiling, one can construct a sequence by
looking at the labels of the tiles. Given a sequence, one can
construct a tiling by associating a tile of length $\ell_i$ to each
letter of type $a_i$, and by concatenating these tiles in the order
given by the sequence. 

More precisely, let $\Xi$ be a space of bi-infinite sequences on a finite alphabet
$\aaa=\{a_1, \ldots, a_n\}$, and let $\ell=(\ell_1, \ell_2, \ldots,
\ell_n)$ be a vector of positive lengths. We then construct a space $\Omega_\ell$
of tilings obtained by associating each letter $a_i$ with a tile
of label $i$ and length $\ell_i$. $\Omega_\ell$ admits  
a natural $\R$ action from translation, while $\Xi$  only admits a
natural $\Z$ action obtained from powers of the shift map.

A natural question is how mixing or spectral properties of the $\Z$
action on $\Xi$ are related to analogous properties of the $\R$ action
on $\Omega_\ell$. When $\ell=(1,\ldots,1)$, this is a simple question,
since the spectrum of translations on $\Omega_\ell$ is simply ($-\sqrt{-1}$
times) the log of the spectrum of the shift operator on
$\Xi$. However, in some tiling spaces changing $\ell$ can
qualitatively change the dynamics, while in others it cannot.

Part of this question was studied in \cite{CS1,CS2}, and is related to
the (\v Cech) cohomology of $\Omega_\ell$. Changes in the lengths of
the tiles are associated with elements of $H^1(\Omega_\ell,
\R)$. There is a subgroup $H^1_{an}(\Omega_\ell,\R)$ of {\em
  asymptotically negligible} length changes, i.e., changes that yield
tiling spaces (and $\R$ actions) that are topologically conjugate to
the original action on $\Omega_\ell$.  Furthermore, a uniform
rescaling of all the tiles, while not a topological conjugacy, does
not change the qualitative mixing and spectral properties. If
$H^1_{an}(\Omega_\ell, \R)$ is a large enough subspace of
$H^1(\Omega_\ell, \R)$, then all length changes are a combination of a
uniform rescaling and a topological conjugacy. This happens, for
instance, for holomogical Pisot substitutions \cite{HPC}.

We call a 1-dimensional tiling space (or the corresponding space of 
sequences) {\em plastic} if all changes to
lengths of the $n$ basic tiles result in topological conjugacies, up
to an overall scale.  We call it {\em totally plastic} if this
remains true for all finite recodings of the tiles. (E.g. collaring,
or dividing a tile into multiple pieces, or amalgamating tiles.) 
Plasticity may be a property of how we label our tiles, but
total plasticity is topological. A tiling space is totally plastic if and
only if $H^1(\Omega_\ell, \R)/H^1_{an}(\Omega_\ell, \R)
= \R$. 

\begin{ex} The space of Fibonacci sequences is generated by the
  substitution $\sigma(a)=ab$, $\sigma(b)=a$. The corresponding space
  of tilings was shown in \cite{RS} to be plastic. To see this,
  suppose that $\Omega_\ell$ and $\Omega_{\ell'}$ are Fibonacci tiling
  spaces corresponding to length vectors $\ell=(\ell_1,\ell_2)$ and
  $\ell' = (\ell'_1, \ell'_2)$. By rescaling $\ell$ or $\ell'$, we can
suppose that $\phi \ell_1 + \ell_2 =
  \phi \ell_1' + \ell_2'$, where $\phi = (1+\sqrt{5})/2$ is the golden
  mean. This implies that $(\ell_1 - \ell_1', \ell_2 - \ell_2')$ is a
  left-eigenvector of the substitution matrix $\left
    ( \begin{smallmatrix} 1 & 1 \cr 1 & 0 \end{smallmatrix} \right )$
  with eigenvalue $1-\phi$, and hence that the difference in length
  between $n$-th order supertiles $\sigma^n(a)$ in $\Omega_\ell$ and
  $\Omega_{\ell'}$ goes to zero as $n \to \infty$, and likewise for
  $\sigma^n(b)$.  Let $\psi_n: \Omega_\ell \to \Omega_{\ell'}$ be a
  map the preserves the underlying sequence, and such that if the
  origin lies a fraction $f$ across an $n$-th order supertile in a
  tiling $T$, then the origin lies the same fraction $f$ across the
  corresponding supertile in $\psi_n(T)$. The tilings $\psi_n(T)$ all
  agree up to translation, and the limiting tiling $\psi(T) := \lim_{n
    \to \infty} \psi_n(T)$ is well-defined. Then $\psi$ is a
  topological conjugacy. (For a more complete description of this
  conjugacy, see \cite{RS} or \cite{TilingsBook}.)
\end{ex}

\subsection{Finitely balanced words}

Let $\Xi$ be a space of sequences. The set of all finite words that
appear within the sequences of $\Xi$ is called the {\em language} of
$\Xi$, and is denoted $\lcal$.

For each letter $a_i$ and each integer $n>0$, let $m_i(n)$ be the
minimum number of times that $a_i$ can appear in a word of length $n$
in $\lcal$, and let $M_i(n)$ be the maximum number of times that $a_i$
can appear.  If there is a constant $C$ such that $M_i(n)-m_i(n) < C$
for all $i$ and all $n$, then we say that the language $\lcal$ is {\em
  finitely balanced}.

\begin{ex} In a Sturmian sequence, the alphabet $\aaa$ consists of two
  letters $a_{1}=a$ and $a_2=b$, and $bb$ is not in the language. We
  then have $m_1(2)=1$ and $M_1(2)=2$, since any word of length 2
  contains either one or two $a$'s. In fact, Sturmian sequences have
  $M_i(n)-m_i(n)=1$ for all $i$ and $n$.
\end{ex}

We can look for variations in the appearance not only of letters, but
of words. Let $w \in \lcal$, and let $m_w(n)$ (resp. $M_w(n)$) be the
minimum (resp. maximum) number of times that $w$ can appear in a word
of length $n$. If for each word $w$ there is a constant $C_w$ such
that $M_w(n)-m_w(n)<C_w$ for all $n$, then we say that $\lcal$ is {\em
  totally finitely balanced}.

\begin{ex} The space of Thue-Morse sequences, coming from the
  substitution $\phi(a)= ab$, $\phi(b) = ba$, is finitely balanced but
  not totally finitely balanced. We have $M_1(n)-m_1(n) \le 2$ for all
  $n$, and likewise $M_2(n)-m_2(n) \le 2$, since all letters of a word
  (excepting possibly the first and/or last letter) pair into
  sub-words $ab$ or $ba$.  However, for words $w$ of length 2,
  $M_{w}(n)-m_{w}(n)$ is unbounded.  On average, 2/3 of the 2-letter
  words in a Thue-Morse sequence are $ab$ or $ba$, while 1/3 are $aa$
  or $bb$, but the language admits words of the form
  $\phi^{2m}(a)\phi^{2m-2}(b) \phi^{2m-4}(a) \cdots$, of length
  $4^m+4^{m-1} + \cdots + 1$, that have $m/3$ too many $ab$'s and
  $ba$'s and $m/3$ too few $aa$'s and $bb$'s.
\end{ex}

Plasticity can allow us to combine insights from sequences and tilings to 
understand the dynamics of both. To date, this has mostly applied to 
substitution sequences (see e.g., \cite{Que, Fogg}) 
and substitution tilings. Substitution tilings have
a natural set of lengths that make the tilings self-similar, and a great 
deal of work, beginning with \cite{Sol}, has been done on the dynamics
of self-similar tilings. Some substitution sequences 
(especially those coming from Pisot substitutions) are known to be
plastic, but others are not. 
Recently, there has been interesting work
(see \cite{Berthe} and references therein) on S-adic sequence spaces. 
This paper arose from an observation that the balancing conditions used in
\cite{Berthe} are exactly the conditions needed for plasticity.

\begin{thm}[Main Theorem] \label{mainthm} A space $\Xi$ of bi-infinite sequences 
  is plastic if and only if its language $\lcal$ 
is finitely balanced. Furthermore, $\Xi$ is
  totally plastic if and only if $\lcal$ is totally finitely
  balanced.
\end{thm}

\section{Subshifts and Tilings}

In this section we review the formalism of 1 dimensional subshifts and
tilings, in preparation for the proof of Theorem \ref{mainthm} in
Section 3. All of the material of this section is standard, and
details may be found in reference books such as \cite{Lind, Fogg,
  TilingsBook}.

Let $\aaa$ be a finite set of symbols, called an {\em alphabet}, and
let $\aaa^*$ be the set of finite words constructed from the letters
of $\aaa$. Then $\aaa^\Z$ is the set of bi-infinite sequences $u =
\cdots u_{-1} u_0 u_1 \cdots$, with each $u_i \in \aaa$. There is a
natural shift action on $\aaa^\Z$, denoted $\sigma$:
\begin{equation} (\sigma(u))_n = u_{n+1}; \qquad (\sigma^{-1}(u))_n =
  u_{n-1}. \end{equation} We give $\aaa^\Z$ a metric where two
sequences are $\epsilon$-close if they agree on $[-1/\epsilon,
1/\epsilon]$. Although this metric gives extra weight to values of a
sequence near the origin, the resulting {\em product topolgy} is
shift-invariant. Equipped with this topology and the $\Z$-action
induced by $\sigma$, $\aaa^\Z$ is called the {\em full shift} on the
alphabet $\aaa$.

A non-empty subset $\Xi \subseteq \aaa^\Z$ is called a {\em subshift}
if it is (1) invariant under $\sigma$ and $\sigma^{-1}$ and (2) closed
in the product topology. The {\em language} $\lcal$ of $\Xi$ is the
set of all finite words that appear within sequences of $\Xi$.

If $u \in \aaa^\Z$, then the {\em orbit closure} of $u$ is the
smallest sub-shift $\Xi_u$ containing $u$. A sequence $v$ is an
element of $\Xi_u$ if and only if every sub-word of $v$ is a sub-word
of $u$. The theory of sub-shifts is often expressed in terms of
specific sequences (typically fixed points of substitutions, as below)
rather than in terms of spaces, since the properties of $\Xi_u$ are
closely related to properties of $u$.

For tilings in 1 dimension, we assume that there are only a finite
number of labels $1, \ldots, n$, and that each label $i$ is associated
with a positive length $\ell_i$. A {\em prototile} is the interval
$[0, \ell_i]$ together with the label $i$. Translating a tile means
translating the interval (say, to $[x, x+\ell_i]$) while preserving
the label.  A {\em tile} is a translated prototile. A {\em tiling} is
a covering of the line by tiles, such that tiles only intersect on
their boundaries.  The translation group $\R$ acts on tilings by
translating all tiles the same distance. 

Given a fixed set of prototiles, we apply a metric to the set of all
tilings such that two tilings $T_1$ and $T_2$ are $\epsilon$-close if
they agree exactly on $[-1/\epsilon, 1/\epsilon]$, up to translations
of each tiling by up to $\epsilon$. A {\em tiling space} is a
non-empty set of tilings that is closed in this topology and invariant
under translation. As with subshifts, the simplest examples of tiling
spaces are closures of translational orbits of individual tilings, 
and we denote the orbit closure of $T$ by $\Omega_T$. 

\subsection{Suspensions}

Let $\Xi$ be a subshift with alphabet $\aaa$.  The {\em suspension} of
$\Xi$ is the mapping cylinder of $\sigma$, namely:
\begin{equation} \Xi \times [0,1] / \sim; \qquad (u,1) \sim
  (\sigma(u),0). \end{equation} This can be viewed as a tiling space
where the set of labels is precisely $\aaa$, and where each tile has
length 1. The point $(u,t)$ consists of a tiling where a tile of label
$u_m$ occupies the interval $[m-t,m+1-t]$. Shifting the sequence is
then equivalent to translating by 1.

If $\ell = (\ell_1, \ldots, \ell_n)$, then the {\em suspension of $\Xi$  by the function} $\ell$ is the quotient of $\Xi \times \R$ by the 
equivalence relation
\begin{equation} (u,t) \sim (\sigma(u), t-\ell_{u_0}) \end{equation}
This is a tiling space in which a tile of type $u_0$ occupies the
interval $[t, t+\ell_{u_0}]$, a tile of type $u_1$ occupies the
interval $[t+\ell_{u_0}, t+\ell_{u_0} + \ell_{u_1}]$, etc. We denote
this space of tilings $\Omega_\ell$.

There is a natural embedding of $\Xi$ in $\Omega_\ell$, namely $u
\mapsto (u,0)$. However, there is no natural map in the reverse
direction, as each translational orbit in $\Omega_\ell$ is connected, 
while $\Xi$ is totally disconnected.

\subsection{Pattern-Equivariant Cohomology}

Let $T$ be a 1-dimensional tiling. We can view this as a decomposition of 
$\R$ into 1-cells (tiles) and 0-cells (vertices). A real-valued 0-cochain is 
an assignment of a real number to each vertex, and a real-valued 1-cochain
is an assignment of a real number to each tile. If $\alpha$ is a 0-cochain,
then $\delta \alpha$ is the coboundary of $\alpha$. That is, if $t$ is a tile
occupying the interval $[a,b]$, then $$\delta \alpha(t) := \alpha (\partial t)
= \alpha(b)-\alpha(a).$$ 
If $\beta = \delta \alpha$, we call $\alpha$ the {\em integral} of $\beta$. 
Any two integrals differ by a constant. 

We say that a cochain is 
{\em strongly pattern-equivariant} (sPE, often abbreviated to PE) 
if the value of a cochain at a cell
depends only on the pattern to a finite radius around that cell. Specifically,
a $k$-cochain $\alpha$ is sPE if there is a radius $R$ with the following 
property:
\begin{itemize}
\item
If two k-cells $c_1$ and $c_2$ are centered at points $x_1$ and $x_2$, 
and if $T-x_1$ and 
$T-x_2$ agree on $[-R,R]$, then $\alpha(c_1)=\alpha(c_2)$. 
\end{itemize}
The first (strongly) PE cohomology of $T$, with real coefficients, is 
\begin{equation} H^1_{PE}(T, \R) := 
\frac{\hbox{1-cochains}}{\partial(\hbox{0-cochains})}.
\end{equation}

The following theorem, originally proved by \cite{KP} for real 
coefficients and extended to arbitrary coefficients in \cite{integer}, 
relates PE cohomology to \v Cech cohomology.

\begin{thm}The PE cohomology of $T$ is 
naturally isomorphic to the \v Cech cohomology of $\Omega_T$.
\end{thm}
 
In particular, if a tiling space is minimal, then the PE cohomology is
the same for all tilings in that space. We then speak of the PE
cohomology of the space, rather than the PE cohomology of any specific
tiling.

A cochain is called {\em weakly pattern equivariant} (wPE) if it is
the uniform limit of sPE cochains. That is, $\alpha$ is wPE if for
each $\epsilon>0$ there exists a radius $R$ such that the value of
$\alpha$ on a cell is determined, to within $\epsilon$, by the pattern
of the tiling within distance $R$ of that cell.

There is a canonical 1-cochain that associates each tile to its length. 
Changes to the lengths of tiles are associated with changes to this 1-cochain,
and hence to an element of $H^1_{PE}$. Conversely, elements of $H^1_{PE}$ can
be represented by sPE 1-cochains. After relabeling tiles according to the 
patterns around them (a process called ``collaring'', which expands 
the alphabet but does not change 
the tiling space as a dynamical system), 
every sufficiently small element of $H^1_{PE}$ can be
associated with a length change. In fact, 

\begin{thm}[\cite{CS2}] Sufficiently small length changes, up to local
equivalence (`MLD'), are parametrized by $H^1_{PE}$. 
\end{thm}

Some of these length changes induce topological conjugacies. The subgroup
of $H^1_{PE}$ that does so is denoted $H^1_{an}$ (for {\em asymptotically negligible}), and is neatly described by a 
theorem of Kellendonk \cite{Kel2}: 

\begin{thm}\label{thm23}
$H^1_{an}$ consists of the classes of sPE 1-cochains
that can be written as the coboundaries of wPE 0-cochains.
\end{thm}

The Gottschalk-Heldlund Theorem \cite{GH} relates bounded integrals to 
continuity in an abstract setting. The following theorem was
proved for tilings in all dimensions in \cite{KS}, 
but in one dimension is in fact an 
immediate corollary of the Gottschalk-Hedlund Theorem.

\begin{thm}\label{thm24}
Let $\beta$ be a sPE 1-cochain. An
integral 
$\alpha$ of $\beta$ is wPE if and only if $\alpha$ is bounded. 
\end{thm}

\section{Proof of Main Theorem}

We first prove the part of the theorem concerning finitely balanced sequences,
and then prove the analogous statements for totally finitely balanced
sequences. 

Suppose that we have a subshift $\Xi$ whose language is finitely balanced. Then
each letter $a_i$ has a well-defined frequency $f_i$, and there is 
a constant $C$ such that 
\begin{equation}\label{lowhi}
 n f_i - C \le m_i(n) \le M_i(n) \le n f_i + C. \end{equation}
Now consider a 1-cochain $\beta_i$ 
on a corresponding tiling (with lengths $\ell=(1,\ldots,1))$ that evaluates
to $1-f_i$ on each tile of type $i$ and to $-f_i$ on each tile that is not 
of type $i$. Note that $\beta_i$ is sPE. 
Pick an integral $\alpha_i$ of $\beta_i$ that equals 0 at a vertex.
Equation \ref{lowhi} implies that $\alpha_i$ takes values between
$-C$ and $C$, and hence by Theorem \ref{thm24} 
that $\alpha_i$ is wPE. Thus, by Theorem \ref{thm23}, changing the lengths of 
the tiles by a multiple of $\beta_i$ induces a topological conjugacy of 
tiling spaces. 

Let $\ell'$ be any vector of positive lengths for the alphabet.
We can write $\ell'$ as a positive multiple $c$ of $\ell$ plus a 
linear combination of the $\beta_i$'s. The length change from $\ell$ to $\ell'$
can then be expressed as an overall scaling by $c$ followed by a sequence 
of topological conjugacies, implying that $\Xi$ is plastic. 

Conversely, suppose that there is a letter $a_i$ 
for which equation (\ref{lowhi}) does
not apply for any choice of $C$ and $f_i$. An arbitrary rescaling,
followed by changing the length $\ell_i$, 
cannot yield a 1-cochain whose integral is bounded. 
Thus changing the length of the $i$-th tile type cannot be 
written as a combination of a topological conjugacy and an overall rescaling.

Next suppose that $\Xi$ is totally finitely balanced. For a given word $w$,
let $\alpha_w$ be an indicator 1-cochain on a corresponding tiling, 
taking the value 1 on a particular 
tile of each occurrence of $w$ and taking the value 0 on all other tiles. 
An equation analogous to (\ref{lowhi}) then applies, with $f_i$ replaced by
the density of the word $w$ and $C$ replaced by $C_w$. Since every sPE 
1-cochain is a finite linear combination of indicator cochains, any sPE
1-cochain whose average value is 0 must have a bounded integral, and so must
have a wPE integral, and so must describe a length change that is a topological
conjugacy. Thus {\em any} length change, after arbitrary collaring, 
can be written as a combination of an overall rescaling and a topological 
conjugacy. 

Finally, suppose that $\Xi$ is not totally finitely balanced, and that
$w$ is a word for which the analogue to (\ref{lowhi}) does not
apply. Collar the tiles so that their labels contain information of a
region of size at least the length of $w$. Consider a length change
wherein the first tile of each occurrence of $w$ is lengthened and all
other tiles are unchanged. That is, the change in $\ell$ is described
by an indicator cochain of $w$. By the failure of (\ref{lowhi}), this
cochain is not a linear combination of the constant cochain $1$ and a
cochain whose integral is wPE (hence bounded), and so the length change is not a
combination of a rescaling and a topological conjugacy.

\noindent {\bf Acknowledgments}  I thank Jarek Kwapisz for pointing out the relation to the Gottschalk-Hedlund Theorem, and all of the participants of the 2014 Workshop on the Pisot Conjecture (at the Lorentz Center in Leiden) for their ideas and encouragement. This work was partially supported by the National Science Foundation under grant  DMS-1101326.


\begin{thebibliography}{ADHM}
\thispagestyle{headings}
\markright{}


%\bibitem[AP]{AP} J.E. Anderson and I.F. Putnam,
%Topological invariants
  %for substitution tilings and their associated $C^*$-algebras,
 %{\em   Ergodic Th.\ \& Dynam.\ Syst.} \textbf{18} (1998), 509--537.

\bibitem[BST]{Berthe} V. Berth\'e, W. Steiner and J. Thuswaldner, 
Geometry, Dynamics and Arithmetic of $S$-adic Shifts, 
preprint 2014 (arXiv:1410.0331)

\bibitem[CS1]{CS1} A.~Clark and L.~Sadun, When Size Matters: Subshifts and 
Their Related Tiling Spaces, {\em Ergodic Theory and Dynamical Systems} {\bf 23}
(2003) 1071--1083.

\bibitem[CS2]{CS2} A.~Clark and L.~Sadun, When Shape Matters: Deformations of Tiling Spaces, 
{\em Ergodic Theory and  Dynamical Systems} {\bf 26} (2006) 69-86.

\bibitem[Fo]{Fogg} N.~P. Fogg, Substitutions in dynamics,
arithmetics and combinatorics. Edited
by V. Berth´e, S. Ferenczi, C. Mauduit and A.
Siegel. Lecture Notes in Mathematics, 1794.
Springer-Verlag, Berlin. (2002)

\bibitem[BBJS]{HPC} M. Barge, H. Bruin,  L. Jones and L. Sadun, 
Homological Pisot Substitutions and Exact Regularity, {\em Israel Journal of Mathematics} {\bf 188} (2012), 281--300.

\bibitem[GH]{GH}  W.H. Gottschalk and G.A. Hedlund. Topological Dynamics. Amer. Math. Soc., Providence,
R. I. (1955).

\bibitem[K1]{KP} J.~Kellendonk, Pattern-equivariant functions and
cohomology, {\em Journal of Physics A.} {\bf 36} (2003), 5765--5772.


\bibitem[K2]{Kel2} J.~Kellendonk, 
 Pattern Equivariant Functions, Deformations and Equivalence of Tiling Spaces, 
 {\em Ergodic Theory and Dynamical Systems}  {\bf 28}  (2008),  no. 4, 1153--1176.

 \bibitem[KS]{KS}
J.~Kellendonk and L.~Sadun, Meyer sets, topological eigenvalues, and Cantor fiber bundles,
{\em Journal of the London Mathematical Society} {\bf 89} (2014) 114--130. 

\bibitem[LM]{Lind} D. Lind and B. Marcus, An Introduction to Symbolic Dynamics and Coding, Cambridge University Press, Cambridge (UK), 1995.

\bibitem[Que]{Que} M. Queff\'elec, Substitution dynamical
systems --- spectral analysis. Lecture Notes
in Mathematics, 1294. Springer-Verlag, Berlin. (1987)

\bibitem[RS]{RS} C.~Radin and L.~Sadun, Isomorphism of Hierarchical Structures,
{\em Ergodic Theory and Dynamical Systems} {\bf 21} (2001) 1239--1248.

\bibitem[Sad1]{integer} L. Sadun,  Pattern-Equivariant Cohomology with Integer Coefficients,
 {\em Ergodic Theory and Dynamical Systems} {\bf 27} (2007), 1991--1998. 

\bibitem[Sad2]{TilingsBook} L. Sadun, Topology of Tiling Spaces. American Mathematical Society, Providence RI, 2008. 

\bibitem[Sol]{Sol} B. Solomyak, Dynamics of Self-Similar Tilings, 
 {\em Ergodic Theory and Dynamical Systems}  {\bf 17}  (1997),  695--738.

\end{thebibliography}
\end{document}